\documentclass{article}
\usepackage{a4,amsmath,amssymb,amsthm,amsfonts}
\usepackage[francais]{babel}
\usepackage[T1]{fontenc}

\begin{document}
\bibliographystyle{plain}

\font\got=eufm10
\font\teneufm=eufm10
\font\seveneufm=eufm7
\font\fiveeufm=eufm5
\newfam\gotfam
\textfont\gotfam=\teneufm
\scriptfont\gotfam=\seveneufm
\scriptscriptfont\gotfam=\fiveeufm
\def\got{\fam\gotfam\teneufm}
\font\tenmsbm=msbm10
\font\sevenmsbm=msbm7
\font\fivemsbm=msbm5
\newfam\spefam
\textfont\spefam=\tenmsbm
\scriptfont\spefam=\sevenmsbm
\scriptscriptfont\spefam=\fivemsbm

%
%
\def\spe{\fam\spefam\tenmsbm}
\def\Bbb{\spe}
\def\Aut{{\rm Aut}}
\def\beginProof{\par{\bf Proof: }}
\def\endProof{${\bf Q.E.D.}$\par}
\def\ar#1{\widehat{#1}}
\def\mn{{\mu_{n}}}
\def\pr{^{\prime}}
\def\prpr{^{\prime\prime}}
\def\mtr#1{\overline{#1}}
\def\ra{\rightarrow}
\def\mQ{{\Bbb Q}}
\def\mR{{\Bbb R}}
\def\mZ{{\Bbb Z}}
\def\mC{{\Bbb C}}
\def\mN{{\Bbb N}}
\def\Qmn{{\mQ(\mn)}}
\def\refeq#1{(\ref{#1})}
\def\umn{^{\mn}}
\def\lmn{_{\mn}}
\def\blb{{\big(}}
\def\brb{{\big)}}
\def\Hom{{\rm Hom}}
\def\Tr{{\rm Tr\,}}
\def\Gal{{\rm Gal}}
\def\End{{\rm End}}
\def\Det{{\rm Det}}
\def\Td{{\rm Td}}
\def\ch{{\rm ch}}
\def\chg{{\rm ch}_{g}}
\def\torus{{\cal T}}
\def\Proj{{\rm Proj}}
\def\Qb{\mtr{\Bbb Q}}
\def\Zn{{{\Bbb Z}/n}}
\def\deg{{\rm deg}}
\def\mod{{\rm mod}}
\def\ac1{\ar{\rm c}_{1}}
\def\NIm#1{{{\rm Im}(#1)}}
\def\NRe#1{{{\rm Re}(#1)}}
\def\rg{{\rm rg}\,}
\def\Id{{\rm Id}}
\def\CX{{\cal X}}
\def\CHL{{\ar{\rm CH}^{*}}}
\def\CH{{\ar{\rm CH}}}
\def\ach{{\ar{\rm ch}}}
\def\Ker{{\rm Ker}}
\def\HDL{H_{\rm Dlb}}
\def\OK{{{\cal O}_K}}
\def\S1{{S}^{1}}
\def\ddc{{\mtr{\partial}\partial\over 2\pi i}}
\def\HG{H_{G}}
\def\chg{{\rm ch}_{g}}
\def\chgt{\widetilde{\rm ch}_{g}}
\def\mathfrakAG{{\widetilde{\mathfrak A}}_{G}}
\def\ZG{Z_{G}}
\def\GS{G^{*}}
\def\GSD{\widehat{\GS}}
\def\Or{{\cal O}}

%
%

%
\def\M#1{\mathbb#1}     
\def\B#1{\bold#1}       
\def\C#1{\mathcal#1}    
\def\E#1{\scr#1}        
\def\R{\M{R}}           %
\def\Z{\M{Z}}           %
\def\N{\M{N}}           
\def\Q{\M{Q}}           %
%
%
\def\ds{\displaystyle}
\def\demo{\noindent {\bf D\'emonstration.}}
\def\hooklongrightarrow{\lhook\joinrel\longrightarrow}  
\def\hooklongleftarrow{\longleftarrow\joinrel\rhook}    
\def\vfi{\varphi}
\def\epsi{\varepsilon}
\def\op#1{\operatorname{#1}}
\def\SP{\op{Spec}\M{Z}}
\def\ssi{si et seulement si}
\def\ov#1{\overline{#1}}
\def\un#1{\underline{#1}}
\def\Res{\op{Res}}
\def\coef{\op{coef}}
\def\div{\op{div}}
\def\Inv{\op{Inv}}
\def\Spec{\op{Spec}}
\def\im{\op{Im}}
\def\rg{\op{rg}}
\def\Hom{\op{Hom}}
\def\Supp{\op{Supp}}
\def\Arg{\op{Arg}}
\def\bi{\binom}

\newtheorem{theor}{Th\'eor\`eme}[section]
\newtheorem{prop}[theor]{Proposition}
\newtheorem{cor}[theor]{Corollaire}
\newtheorem{lemma}[theor]{Lemme}
\newtheorem{sublem}[theor]{sous-lemme}
\newtheorem{defin}[theor]{Definition}
\newtheorem{conj}[theor]{Conjecture}
\theoremstyle{definition}
\newtheorem*{remarque}{Remarque}

\author{Vincent MAILLOT\footnote{Institut de Math\'ematiques de Jussieu,
     Universit\'e Paris 7 Denis Diderot,
     Case Postale 7012,
     2 place Jussieu,
     F-75251 PARIS CEDEX 05,
     {courriel} : {\em vmaillot@math.jussieu.fr}}\ \ et
     Damien ROESSLER\footnote{Institut de Math\'ematiques de Jussieu,
     Universit\'e Paris 7 Denis Diderot,
     Case Postale 7012,
     2 place Jussieu,
     F-75251 PARIS CEDEX 05,
     {courriel} : {\em roessler@math.jussieu.fr}}}
\title{Conjectures sur les d\'eriv\'ees
logarithmiques des fonctions $L$ d'Artin aux entiers n\'egatifs}
\maketitle
\begin{abstract}
Nous formulons plusieurs variantes d'une conjecture reliant
le degr\'e arith\-m\'etique de certains fibr\'es hermitiens
aux valeurs prises aux entiers n\'egatifs par la d\'eriv\'ee logarithmique
des fonctions $L$ d'Artin, g\'en\'eralisant des conjectures
de Colmez et Gross-Deligne et compl\'ementant ainsi les conjectures 
de Beilinson pour les motifs d'Artin. 
Nous annon\c cons plusieurs r\'esultats en
direction de ces \'enonc\'es.
\end{abstract}
\date
\begin{flushleft} {\bf Conjectures on the logarithmic
derivatives of Artin $L$-functions at non-positive integers}
\end{flushleft}
\indent{\footnotesize {\bf Abstract} {\it -- 
We formulate several variants of a conjecture relating
the arithmetic degree of certain hermitian fibre bundles with
the values of the logarithmic derivative of Artin's $L$-functions at 
negative integers. This generalizes conjectures by 
Colmez and Gross-Deligne and complements Beilinson's conjectures
for the Artin motives. We announce several results 
in the direction of these statements. 
}}
\date
\parindent=0pt
\section{Pr\'eliminaires}
On appelle vari\'et\'e arithm\'etique tout sch\'ema r\'egulier
qui est quasi-projectif et plat sur un anneau arithm\'etique
(au sens de Gillet-Soul\'e \cite{GS2}). Soient $R$ un tel anneau
et $\cal C$ l'ensemble
des couples $(X,D)$, o\`u $f:X\ra\Spec R$ est une vari\'et\'e arithm\'etique
sur $R$ et $D$ un diviseur de Weil \`a croisements normaux de $X(\mC)$. Un {\it fibr\'e
hermitien sur $X$ \`a singularit\'es logarithmiques le long de $D$}
est la donn\'ee d'un fibr\'e coh\'erent $E$ sur $X$ localement libre
sur $X(\mC)$ et d'une m\'etrique hermitienne sur
$E_\mC|_{X(\mC)\setminus D}$ qui est {\it bonne}
le long de $D$ au sens de \cite{Mu} ; on note $\overline{E}$ un tel couple.
Dans tout ce qui suit
$\ar{\rm CH}^*(X)$ d\'esigne l'anneau de Chow arithm\'etique usuel
associ\'e \`a $X$ (cf. \cite{GS2}) et $(X,D)\mapsto\CHL(X,D)$
la fl\`eche construite par J. I. Burgos et U. K\"uhn
de $\cal C$ vers la cat\'egorie des groupes ab\'eliens $\mN$-gradu\'es 
(\`a para\^\i tre). 
Cette th\'eorie est une extension en dimension 
sup\'erieure des th\'eories de Bost (voir \cite{Bo}) et K\"uhn 
(voir \cite{Ku}) pour les surfaces. On obtient cette 
fl\`eche en rempla\c cant dans une construction 
de Burgos (voir \cite{Bu2}) des formes lisses par des formes 
\`a singularit\'es faibles le long de $D$.
Elle v\'erifie (entre autres choses) les propri\'et\'es
suivantes~:
\smallskip
\begin{description}
\item[(1)] le groupe $\CHL(X,D)\otimes_\Z\Q$ est muni d'une structure 
naturelle d'anneau commutatif unitaire pour laquelle la $\N$-graduation 
induite est une graduation d'anneau; 
\item[(2)] si $X$ est propre
sur $\Spec R$, il existe une application \og image directe\fg\ $f_{*}:
\ar{\rm CH}^{d_{X}+1}(X,D)\ra\ar{\rm CH}^{1}(R)$, o\`u
$d_{X}$ est la dimension relative de $X$ sur $R$; cette application 
est un homomorphisme de groupes;
\item[(3)] \`a tout entier $r\geqslant 0$ et tout
fibr\'e hermitien $\mtr{E}$ sur $X$
\`a singularit\'es logarithmiques le long de $D$, on
 peut associer une $r$-i\`eme \og classe de Chern\fg\
$\ar{c}_{r}(\mtr{E})\in\ar{\rm CH}^r(X,D)$;
\item[(4)] \`a tout $R$-morphisme de sch\'emas
$g:X\pr\ra X$ transverse \`a $D$, o\`u $X\pr$ est une vari\'et\'e arithm\'etique
sur $R$, on peut associer un morphisme de groupes gradu\'es \og image
r\'e\-ci\-pro\-que\fg\
$g^*:\CHL(X,D)\ra \ar{\rm CH}^{*}(X\pr,g^{-1}D)$ qui est 
un morphisme d'anneaux lorsqu'on tensorise par $\mQ$; 
\item[(5)] l'\'egalit\'e $g^{*}(\ar{c}_{r}(\mtr{E}))=
\ar{c}_{r}(g^{*}(\mtr{E}))$ vaut pour tout $r\geqslant 0$;
\item[(6)] on dispose d'un morphisme \og d'oubli\fg\
$\zeta:\CHL(X,D)\ra{\rm CH}^*(X)$ compatible
\`a l'application image
r\'e\-ci\-pro\-que du (4) et \`a celle usuelle de la
th\'eorie de Chow, ainsi qu'\`a la graduation et \`a la formation des classes
de Chern;
\item[(7)] on dispose pour tout entier $p\geqslant 0$ d'un morphisme
$\omega:\CH^{p-1}(X,D)\ra Z^{p,p}(X(\mC),D)$
et d'un complexe de groupes~:
$$
\bigoplus_{p\geqslant 0}H_{\cal D}^{2p-1}(X_\mR,\mR(p))\stackrel{a}{\ra}
\CHL(X,D)\otimes_\mZ\mQ\stackrel{(\zeta\oplus\omega)}{\ra}
{\rm CH}^{*}(X)\otimes_\mZ\mQ\oplus
\bigoplus_{p\geqslant 0}Z^{p,p}(X({\mC}),D)
$$
o\`u $H_{\cal D}^{2p-1}(X_\mR,\mR(p))$ d\'esigne la partie
invariante par conjugaison complexe de
la cohomologie de Deligne r\'eelle $H_{\cal D}^{2p-1}(X(\mC),\mR(p))$ et
$Z^{p,p}(X(\mC),D)$ le $\mC$-espace vectoriel des formes r\'eelles
de type $(p,p)$ invariantes par conjugaison et {\it bonnes}
le long de $D$ au sens de \cite{Mu};
\item[(8)] si l'on munit $\bigoplus_{p\geqslant 0}Z^{p,p}(X({\mC}),D)$
de la structure d'anneau d\'efinie par le produit des formes,
alors le morphisme $\zeta\oplus\omega$ est
un morphisme d'anneaux ; l'image par $a$ de
$\bigoplus_{p\geqslant 0}H_{\cal D}^{2p-1}(X_\mR,\mR(p))$
est un id\'eal de carr\'e nul ; et si l'on note $c:{\rm CH}^{*}(X)\ra
\bigoplus_{p\geqslant 0}H_{\cal D}^{2p}(X_\mR,\mR(p))$ l'application
cycle, alors $a(x)\cdot y=a(c(\zeta(y))\cdot x)$, o\`u
$x\in \bigoplus_{p\geqslant 0}H_{\cal D}^{2p-1}(X_\mR,\mR(p))$,
$y\in\CHL(X,D)\otimes_\mZ\mQ$ et le premier point $\cdot$ (resp. le
deuxi\`eme) dans la derni\`ere \'egalit\'e d\'esigne le produit
dans l'anneau $\CHL(X,D)\otimes_\mZ\mQ$ (resp. en cohomologie de Deligne);
\item[(9)] enfin si $D$ est vide et $X$ projective, il existe un isomorphisme canonique
de groupes gradu\'es $\CHL(X,D)\simeq\ar{\rm CH}^*(X)$ qui 
est un isomorphisme d'anneaux lorsqu'on tensorise par $\mQ$ et
via cet isomorphisme les classes $\ar{c}_{r}(\cdot)$
et les morphismes $f_{*}$, $g^{*}$, $\zeta$, $a$ et $\omega$
introduits ci-dessus
co\"\i ncident avec les objets \'eponymes de la th\'eorie
de Gillet-Soul\'e (cf. \cite{GS2} et \cite{GS3}).
\end{description}
Si plus g\'en\'eralement on consid\`ere dans ce qui pr\'ec\`ede des cycles \`a 
coefficients dans un corps $A\subseteq\mR$ plut\^ot 
que dans $\mZ$, on obtient une fl\`eche
vers la cat\'egorie des $A$-alg\`ebres $\mN$-gradu\'ees que 
l'on note alors $(X,D)\mapsto\CHL(X,D)_A$. On veillera 
\`a ne pas confondre $\CHL(X,D)_A$ avec $\CHL(X,D) \otimes_\mZ A$.

\section{Une conjecture pour les fibrations semi-ab\'eliennes}

Soient $p:{\cal A}\ra B$ un sch\'ema semi-ab\'elien lisse  
de dimension relative $d_{\cal A}$ au-dessus de $B$ une vari\'et\'e 
arithm\'etique de dimension $d$ sur $R$ et $L$ un fibr\'e 
en droites sur $\cal A$ ample relativement \`a $p$. 
On suppose qu'il existe 
un ouvert non vide $U\subseteq B$ au-dessus duquel la restriction de
$\cal A$ est un sch\'ema ab\'elien et, par abus de langage, on 
note encore $U$ le plus grand ouvert ayant cette propri\'et\'e. 
On note $S$ le compl\'ementaire de $U$ dans $B$ dont on suppose
que $D:=S(\mC)$ est un diviseur \`a croisements normaux dans $B(\mC)$ et  
l'on fait l'hypoth\`ese que le sch\'ema ab\'elien dual de ${\cal A}_{|U}$ s'\'etend 
en un sch\'ema semi-ab\'elien sur $B$ que l'on note ${\cal A}^\vee$.
Le morphisme $f_L:{\cal A}_{|U}\ra{{\cal A}}^\vee_{|U}$ induit par $L$
s'\'etend en un morphisme de ${\cal A}$ vers ${\cal A}^\vee$ que 
l'on note encore $f_L$.

Soit $K$ une extension finie de $\mQ$ dont on note 
${\cal O}_K$ l'anneau des entiers et pour lequel
on suppose 
que chacun des plongements d'anneaux de 
$\OK$ dans $\mC$ se factorise par  $R$. 
On suppose donn\'e un plongement 
  d'anneaux $\iota:{\cal O}_{K}
\hookrightarrow \End_{B}({\cal A})$ o\`u $\End_{B}({\cal A})$ d\'esigne
l'anneau des endomorphismes du sch\'ema en groupes $p:{\cal A}\ra B$. Par abus 
de langage, on \'ecrit  
$\Omega_{\cal A}$ (resp. $\Omega_{{\cal A}^\vee}$) pour l'image r\'eciproque par 
la section unit\'e du faisceau des diff\'erentielles de 
$\cal A$ (resp. ${\cal A}^\vee$) sur $B$. 
Le fibr\'e $\Omega_{{\cal A},\mC}$ 
est muni de la m\'etrique $L^2$ induite 
par $L_\mC$ et 
$\Omega_{{\cal A}^\vee,\mC}$ h\'erite de cette m\'etrique via 
l'isomorphisme $f_{L,\mC}^{*}:\Omega_{{\cal A}^\vee,\mC}|_{U_\mC}\ra  
\Omega_{{\cal A},\mC}|_{U_\mC}$ induit par $f_{L}$. 
La m\'etrique canonique sur $\Omega_{{\cal A},\mC}$ (resp. sur
$\Omega_{{\cal A}^\vee,\mC}$) est bonne le long de $D$ 
(voir \cite{Mu}). 
Soit $\sigma\in\Hom(\OK,R)$ et $b\in{\cal O}_{K}$; on note 
$\Omega_{{\cal A},\sigma,b}$ (resp. $\Omega_{{\cal A}^\vee,\sigma,b}$) 
le noyau du morphisme de faisceaux  
$\iota(b)-\sigma(b):\Omega_{\cal A}\ra\Omega_{\cal A}$ (resp. 
$\iota(b)-\sigma(b):
\Omega_{{\cal A}^\vee}\ra\Omega_{{\cal A}^\vee}$) et l'on note 
$\Omega_{{\cal A},\sigma}$ (resp. $\Omega_{{\cal A}^\vee,\sigma}$) le faisceau 
$\cap_{b\in{\cal O}_{K}}\Omega_{{\cal A},\sigma,b}$ 
(resp. $\cap_{b\in{\cal O}_{K}}
\Omega_{{\cal A}^\vee,\sigma,b}$). Ce dernier 
est un faisceau coh\'erent car 
${\cal O}_{K}$ est finiment engendr\'e comme 
$\mZ$-module. On suppose que $\Omega_{{\cal A},\sigma}$ (resp. 
$\Omega_{{\cal A}^\vee,\sigma}$) est localement 
libre sur $B$, et que les fibr\'es $\Omega_{{\cal A},\sigma,\mC}$ 
(pour tout $\sigma \in\Hom(\OK,R)$) 
sont orthogonaux entre eux fibre \`a fibre (il en est alors 
de m\^eme des $\Omega_{{\cal A}^\vee,\sigma,\mC}$).
La m\'etrique canonique sur $\Omega_{{\cal A},\mC}$ 
(resp. de $\Omega_{{\cal A}^\vee,\mC}$) induit par restriction sur les
$\Omega_{{\cal A},\sigma,\mC}$ 
(resp. les $\Omega_{{\cal A}^\vee,\sigma,\mC}$) une m\'etrique
qui d'apr\`es ce qui pr\'ec\`ede est bonne le long de $D$.
Enfin on note ${\cal H}_{n}:=\sum_{j=1}^{n}{1/j}$ le 
$n$-i\`eme nombre harmonique et, pour tout caract\`ere d'Artin $\chi$ de $K$,
on note $\mQ(\chi)$ le corps engendr\'e par l'ensemble 
des valeurs de $\chi$ et $\mQ^+(\chi)$ le sous-corps 
$\mQ(\chi)\cap\mR$.
\begin{conj}  
Pour tout $n\geqslant 1$, l'\'egalit\'e~:
\begin{eqnarray*}
\lefteqn{
\frac{1}{2}\sum_{\sigma\in\Hom(\OK,R)}\ach^{[n]}(\mtr{\Omega}_{{\cal A},\sigma}\oplus 
\mtr{\Omega}_{{\cal A}^\vee,\sigma}^{\vee})
\chi(\sigma)}\\&=&
-a\blb\blb {L\pr(\chi,1-n)\over L(\chi,1-n)}+{1\over 2}{\cal H}_{n-1}
-{c_\chi \log(2)\over  1-2^{-n}}\brb
\sum_{\sigma\in\Hom(\OK,R)}\ch^{[n-1]}
(\Omega_{{\cal A},\sigma,\mC})\brb\chi(\sigma)
\end{eqnarray*}
vaut dans $\CH^{n}(B,D)_{\mQ^+(\chi)}\otimes_{\mQ^+(\chi)}\mQ(\chi)$  
pour tout caract\`ere d'Artin irr\'eductible 
$\chi$ de $K$ tel que $L(\chi,1-n)\not= 0$, sauf \'eventuellement lorsque 
$n=1$ et $\chi$ est le caract\`ere 
trivial. Ici $c_{\chi}=1$ si $\chi$ est le
caract\`ere trivial  et $c_\chi=0$ sinon.
\label{conjI}
\end{conj}

{\bf Remarques} 

(1) L'\'enonc\'e que l'on obtient en appliquant 
le morphisme d'oubli $\zeta$ aux deux c\^ot\'es de la conjecture \ref{conjI} 
est (\`a la connaissance des auteurs) \'egalement conjectural. 

(2) Supposons que $R$ est un ouvert d'un anneau de corps de
nombres et soit $\Spec \mZ[1/N]$ le plus 
grand ouvert de $\mZ$ au-dessus duquel $\Spec R$ est fini ; lorsque $B$ est propre,
on \'ecrit $\ar{\deg}:\ar{\rm CH}^{d+1}(B,D)_{\mQ^+(\chi)}
\ra\ar{\rm CH}^{1}(\Spec \mZ[1/N])_{\mQ^+(\chi)}$ 
pour la compos\'ee des applications \og image directe\fg\ 
$\ar{\rm CH}^{d+1}(B,D)_{\mQ^+(\chi)}\ra\ar{CH}^{1}(R)_{\mQ^+(\chi)}$ et 
$\ar{CH}^{1}(R)_{\mQ^+(\chi)}\ra
\ar{\rm CH}^{1}(\Spec \mZ[1/N])_{\mQ^+(\chi)}$ introduites au 1.
De l'\'egalit\'e obtenue en appliquant 
le morphisme 
$\ar{\deg}$ aux deux c\^ot\'es de la formule conjectur\'ee 
pour $n = d + 1$, on d\'eduit une \'egalit\'e de nombres r\'eels \`a 
des sommes de termes de la forme $b_\chi \log q$ pr\`es, 
o\`u $q$ est un entier
positif divisant $N$ et $b_\chi$ un \'el\'ement de $\mQ(\chi)$.
En particulier, si 
$R$ est un anneau de corps de nombres, on obtient une \'egalit\'e 
de nombres r\'eels.
\vskip 0.1cm

Lorsque $K$ est un corps CM et que $n = 1$, 
$d_{\cal A}=[K:\mQ]/2$ et $S$ est vide, la conjecture \ref{conjI}
est un l\'eger raffinement de la conjecture de Colmez \cite[Conj. 0.4, p. 632]{C}
qu'il prouve lorsque $K$ est une extension 
ab\'elienne de $\mQ$ (\`a un facteur $\log(2)$ pr\`es). Le cas 
de \ref{conjI} o\`u $K$ est un corps 
quadratique imaginaire, $n=1$ et $S$ est vide implique 
un raffinement d'un th\'eor\`eme 
de Gross \cite[Th. 3, Par. 3, p. 204]{Gross}. 
Si $2$ est inversible dans $R$, le calcul \cite[Th. 6.1, p. 229]{Ku} de K\"uhn  
(voir aussi celui fait ind\'ependamment par Bost) 
\'etablit la conjecture \ref{conjI} pour 
$K=\mQ$, $n=2$ et $d_{\cal A}=1$ et sa 
conjecture \cite[Cor. 2.9.14]{KuThe} 
implique la conjecture \ref{conjI} lorsque $K$ est un corps 
quadratique totalement r\'eel, $n=2$ et 
$d_{\cal A}=2$ (voir la proposition \ref{corcp}). Le cas de \ref{conjI} 
o\`u $K$ est un corps quadratique imaginaire, $d_{A}=2$ 
et $R$ est un certain anneau arithm\'etique peut \^etre 
d\'eduit d'une conjecture de Kudla, Rapoport et Yang sur les d\'eriv\'ees 
de certaines s\'eries d'Eisenstein incoh\'erentes (voir \cite[(0.16)]{KRY}) ;
et les cas o\`u $K$ est un corps 
quadratique totalement r\'eel, $d_{\cal A}=2$ et 
$n=2$ (voir la proposition \ref{corcp})
et o\`u $K = \mQ$, $d_{A}=2$, $d =3$ et $n=2$, $4$
permettent de pr\'eciser une conjecture de Kudla (voir \cite[(6.11) et (6.12)]{Kudla}).
Enfin le cas o\`u  
$K=\mQ$ et $S$ est vide est conjectur\'e par K\"ohler dans 
\cite[Th. 3.3]{K1}. La proposition suivante donne 
comme application de ce qui pr\'ec\`ede
une formule pour le degr\'e arithm\'etique des fibr\'es 
de formes modulaires sur les 
vari\'et\'es modulaires de Hilbert, munis de leur m\'etrique de Petersson. 
\begin{prop}
Soit ${\cal A}\ra B$ une fibration semi-ab\'elienne lisse munie d'une action
$\iota:\OK\hookrightarrow\End_{B}({\cal A})$
de l'anneau des entiers d'un corps de nombres ab\'elien totalement r\'eel $K$ 
et satisfaisant aux hypoth\`eses 
de la conjecture \ref{conjI}. Supposons de plus que $B$ est propre,
$d_{\cal A} = [K:\mQ] = d$ et que 
$R$ est un ouvert d'un anneau de 
corps de nombres o\`u le degr\'e 
de $f_L$ est inversible. Soit 
$\mZ[1/N]$ le plus grand ouvert au-dessus duquel 
$R$ est fini et soit $K_A$ le corps engendr\'e par les valeurs 
de tous les caract\`eres d'Artin (pairs) de $K$. 
L'\'egalit\'e~:
$$
\ar{\deg}(\ar{c}_{1}^{d+1}(\mtr{\Omega}_{\cal A}))=
-(d+1)\Bigl( {d\over 3}
{\zeta\pr_\mQ(-1)\over \zeta_\mQ(-1)}+{2\over 3}
{\zeta_{K}\pr(-1)\over \zeta_{K}(-1)}+{d\over 2}-
{4(d+2) \over 9}\log(2)\Bigr)\,\deg(\Omega_{{\cal A},\mC})
$$ 
dans $\ar{\rm CH}^{1}(\Spec\mZ[1/N])_{K_A^+}$ 
est une cons\'equence de la conjecture \ref{conjI} pour 
$n=2$ et $\chi$ parcourant les caract\`eres d'Artin (pairs)  
de $K$. 
\label{corcp}
\end{prop}
\begin{remarque}
Il ne semble pas possible de d\'eduire la formule pr\'ec\'edente
de la conjecture \ref{conjI} lorsque $K$ n'est plus 
une extension ab\'elienne de $\mQ$.
On peut n\'eanmoins se demander si la proposition \ref{corcp} reste 
vraie sans cette hypoth\`ese.
\end{remarque}
La proposition suivante donne 
comme application de ce qui pr\'ec\`ede
une formule pour le degr\'e arithm\'etique des fibr\'es 
de formes modulaires sur les 
vari\'et\'es modulaires classifiant les surfaces 
ab\'eliennes \`a multiplication quaternionique, 
munis de leur m\'etrique de Petersson.  
\begin{prop}
Soit ${\cal A}\ra B$ une fibration semi-ab\'elienne lisse   
munie d'une action
$\iota:\OK\hookrightarrow\End_{B}({\cal A})$ 
d'un anneau de corps de nombres quadratique imaginaire 
 et satisfaisant aux hypoth\`eses 
de la conjecture \ref{conjI}. Supposons de plus que $B$ est propre,
$d_{\cal A} = [K:\mQ] = d+1 = 2$ et 
que 
$R$ est un ouvert d'un anneau de 
corps de nombres o\`u le degr\'e 
de $f_L$ est inversible. 
Soit 
$\mZ[1/N]$ le plus grand ouvert au-dessus duquel 
$R$ est fini. 
L'\'egalit\'e~:
$$
\ar{\deg}(\ar{c}_{1}^{2}(\mtr{\Omega}_{\cal A}))=
-(4{\zeta\pr_\mQ(-1)\over\zeta_\mQ(-1)}-{16\over 3}\log(2)+2)\,\deg(\Omega_{{\cal A},\mC})
$$
dans $\ar{\rm CH}^{1}(\Spec\mZ[1/N])_\mQ$ 
est une cons\'equence de la conjecture \ref{conjI} pour $n=2$ et 
$\chi$ le caract\`ere trivial de $K$ et pour $n=1$ et $\chi$ l'unique 
caract\`ere d'Artin non-trivial de $K$. 
\end{prop}
A la fin d'une nouvelle version de \cite{K1}, K\"ohler conjecture une formule 
pour le degr\'e arithm\'etique des fibr\'es 
de formes modulaires sur les 
vari\'et\'es modulaires de Siegel. Cette formule est une 
cons\'equence de la conjecture \ref{conjI}, pour $n$ parcourant 
les entiers pairs positifs et $\chi$ le caract\`ere 
trivial.\\
Indépendamment de ce qui précède, J. Kramer a conjecturé (conversations 
avec le second auteur) que le degr\'e 
arithm\'etique du fibr\'e des formes modulaires sur 
une vari\'et\'e de Shimura est donn\'e par une combinaison 
lin\'eaire rationnelle de d\'eriv\'ees logarithmiques 
de fonctions $\zeta$ de corps de nombres \'evalu\'ees 
en des entiers n\'egatifs. Les deux
propositions précédentes sont compatibles avec
une telle conjecture.  

\section{Une conjecture pour les motifs relatifs}

Soient $Y$ une vari\'et\'e arithm\'etique lisse    
sur $R$ un anneau arithm\'etique de 
Dedekind et $\cal V$ un motif relatif sur $Y$ \`a coefficients 
dans $\mZ[{1\over m}]$ (avec $m$ un entier positif) et effectif pour 
l'\'equivalence rationnelle. Le motif $\cal V$ consiste par 
d\'efinition en un couple form\'e d'une fibration projective et lisse $f:X\ra Y$ et 
d'un idempotent pour la composition des correspondances 
dans ${\rm CH}^*(X\times_Y X,\mZ[{1\over m}])$; nous 
renvoyons \`a \cite{DeMu} 
pour la d\'efinition g\'en\'erale de la cat\'egorie des motifs relatifs sur $Y$.  
On se donne sur $X$ un fibr\'e en droites $L$ ample 
relativement \`a $f$. Nous \'ecrirons 
$\HDL({\cal V})$ pour d\'esigner la cohomologie 
de Dolbeaut de $\cal V$ (on rappelle que 
si l'idempotent associ\'e \`a $\cal V$ est l'identit\'e, 
$\HDL({\cal V})$ est le faisceau $\oplus_{k\geq 0}\oplus_{p+q=k}
R^q f_*\Lambda^p\Omega(f)$); c'est un faisceau 
coh\'erent $\mN$-gradu\'e sur $Y$. On suppose que le faisceau 
$\HDL({\cal V})_\mC$ est localement libre; il est alors 
muni de la m\'etrique $L^2$ induite par $c_{1}(L_\mC)\in 
\HDL^2(X/Y)_\mC$.\\
Soit $P(T)$ un polyn\^ome unitaire \`a coefficients dans $\mZ[{1\over m}]$ et 
irr\'eductible sur $\mQ$. On pose
$\Or:=\mZ[{1\over m}][T]/(P(T))$ et l'on note $K$ le corps 
des fractions de $\Or$ (c'est une extension finie de 
$\mQ$) et l'on suppose que chacun des plongements d'anneaux de 
$\Or$ dans $\mC$ se factorise par $R$.
On suppose donn\'e un morphisme 
d'anneaux $\iota:\mZ[{1\over m}][T]\ra\End({\cal V})$
dont le noyau est contenu dans l'id\'eal $(P(T))$.
Pour tout $\sigma\in\Hom(\Or,R)$ on note 
$\HDL({\cal V})_{\sigma}$  
le noyau du morphisme de faisceaux  
$\iota(T)-\sigma(T):\HDL({\cal V})\ra\HDL({\cal V})$ (qui est coh\'erent); 
on suppose que 
$\HDL({\cal V})_\sigma$ est localement libre 
et l'on munit $\HDL({\cal V})_{\sigma,\mC}$ de la m\'etrique 
induite par restriction de celle de $\HDL({\cal V})_\mC$. On suppose 
de plus que les fibr\'es $\HDL({\cal V})_{\sigma,\mC}$ (pour tout $\sigma\in\Hom(\Or,R)$) 
sont orthogonaux entre eux fibre \`a fibre. 
Enfin on note $r$ le produit de tous les nombres premiers $q$ au-dessus desquels
il existe un id\'eal premier $\cal Q$ de $\Or_K$ 
tel que $\Or_{\cal Q}$ ne co\"\i{}ncide pas avec $\Or_{K,{\cal Q}}$. 
\begin{conj}
Pour tout $n\geqslant 1$ et $k\geq 0$, l'\'egalit\'e~:
\begin{eqnarray*}
\lefteqn{
{1\over 2}\sum_{\sigma\in\Hom(\Or,R)}
\ach^{[n]}({{ H}^{k}_{\rm Dlb}(\overline{\cal V})}_\sigma)
\chi(\sigma)}\\&=&
-\sum_{\sigma\in\Hom(\Or,R)}a\blb( {L\pr(\chi,1-n)\over L(\chi,1-n)}+{1\over 2}{\cal H}_{n-1}-
{c_\chi \log(2)\over  1-2^{-n}})
\sum_{p+q=k}p\cdot \ch^{[n-1]}
(\HDL^{p,q}({\cal V})_{\sigma, \mC})\brb\chi(\sigma)\\
\end{eqnarray*}
vaut dans $\CH^{n}(Y[{1\over mr}])_{\mQ^+(\chi)}
\otimes_{\mQ^+(\chi)}\mQ(\chi)$ 
pour tout $k \geqslant 0$ et tout caract\`ere d'Artin irr\'eductible 
$\chi$ de $K$ tel que $L(\chi,1-n)\not=0$, 
sauf \'eventuellement lorsque $n=1$
et $\chi$ est le caract\`ere  trivial. Ici $c_{\chi}=1$ si $\chi$ est le
caract\`ere trivial  et $c_\chi=0$ sinon. 
\label{conjII}
\end{conj}
{\bf Remarques}

(1) La conjecture \ref{conjII} implique 
la conjecture \ref{conjI} lorsque $Y=\Spec R$ 
et que $R$ est un corps. 

(2) L'\'enonc\'e g\'eom\'etrique
obtenu en appliquant le morphisme d'oubli 
$\zeta$ aux deux c\^ot\'es 
de la conjecture est (\`a la connaissance 
des auteurs) \'egalement conjectural. 
On peut interpr\'eter cet \'enonc\'e comme 
une version relative de la formule des traces de Lefschetz pour la 
cohomologie singuli\`ere \`a coefficients dans $\mC$. 

(3) Lorsque 
$Y=\Spec \Qb$, la conjecture \ref{conjII} est une extension non-ab\'elienne
d'une variation de la conjecture 
\og des p\'eriodes\fg\ de Gross-Deligne \cite[p. 205]{Gross}.
 
On remarquera \'egalement 
qu'au vu de la conjecture \ref{conjI}, on souhaiterait 
\'etendre la conjecture \ref{conjII} \`a des motifs associ\'es 
\`a des fibrations 
\`a singularit\'es semi-stables ou, en d'autres termes, \`a une 
certaine sous-cat\'egorie de la \og cat\'egorie des faisceaux en motifs mixtes\fg .\\
La conjecture suivante est une version affaiblie de la 
conjecture \ref{conjII} dont elle est une cons\'equence 
imm\'ediate~:
\begin{conj} 
Pour tout $n\geqslant 1$, l'\'egalit\'e~:
\begin{eqnarray*}
\lefteqn{
{1\over 2}\sum_{k\geqslant 0}(-1)^{k}\sum_{\sigma\in\Hom(\Or,R)}
\ach^{[n]}({{H}^{k}_{\rm Dlb}(\overline{\cal V})}_\sigma)
\chi(\sigma)}\\&=&
-\sum_{\sigma\in\Hom(\Or,R)}a\blb\sum_{k\geqslant 0}(-1)^{k}\blb {L\pr(\chi,1-n)\over 
L(\chi,1-n)}+
{1\over 2}{\cal H}_{n-1}-
{c_\chi \log(2)\over  1-2^{-n}}\brb
\sum_{p+q=k}p\cdot\ch^{[n-1]}
(\HDL^{p,q}({\cal V})_{\sigma, \mC})\brb\chi(\sigma)\\
\end{eqnarray*}
vaut dans $\CH^{n}(Y[{1\over mr
}])_{\mQ^+(\chi)}\otimes_{\mQ^+(\chi)}\mQ(\chi)$ 
pour tout caract\`ere d'Artin irr\'eductible 
$\chi$ de $K$ tel que 
$L(\chi,1-n)\not=0$, 
sauf \'eventuellement lorsque $n=1$ et $\chi$ est le
caract\`ere  trivial. Ici $c_{\chi}=1$ 
si $\chi$ est le caract\`ere trivial 
et $c_\chi=0$ sinon.
\label{conjIIa}
\end{conj}

\section{R\'esultats}
 
Les hypoth\`eses et les notations utilis\'ees ci-dessous sont celles 
du paragraphe pr\'ec\'e\-dent.
Les d\'emonstrations des r\'esultats annonc\'es ici seront 
publi\'ees ult\'erieurement.
\vskip 0.1cm

Dans tout ce qui suit, si $\chi$ est un caract\`ere d'Artin et 
$H$ un corps de nombres contenant $\mQ(\chi)$, 
nous dirons que la conjecture \ref{conjI} (resp. \ref{conjII}, 
resp. \ref{conjIIa}) vaut pour $\chi$ {\it apr\`es produit tensoriel par 
$H$} si la conjecture obtenue en rempla\c cant $\mQ(\chi)$ par $H$ dans 
l'\'enonc\'e correspondant est v\'erifi\'ee.\\
Le théorème suivant relie la conjecture \ref{conjIIa}
\`a une formule de Lefschetz relative conjectur\'ee dans \cite[Appendix]{KR1}
et maintenant démontrée grâce à \cite{KR1} et \cite{BiMa}. 
\begin{theor}
Supposons que $\cal V$ est le motif relatif associ\'e 
\`a une fibration lisse $f:X\ra Y$, que 
$P(T)$ est le $q$-i\`eme polyn\^ome cyclotomique et
qu'il existe une racine primitive $q$-i\`eme de l'unit\'e dans 
$K = \mQ(\mu_q)$ dont l'image par $\iota$ est repr\'esent\'ee  
par un automorphisme d'ordre $q$ de la fibration $f$. Supposons 
de plus que $q$ est inversible dans $R$ et, pour tout $\chi$ caract\`ere de 
Dirichlet primitif modulo $q$, notons
$\mQ(\chi,\mu_q)$ le corps engendr\'e par les valeurs
de $\chi$ et les racines $q$-i\`emes de l'unit\'e.
La conjecture \cite[Appendix]{KR1} implique que la 
conjecture \ref{conjIIa} vaut pour $\chi$
apr\`es produit tensoriel par $\mQ(\chi,\mu_q)$.  
\label{res1}
\end{theor}
{\bf Esquisse de la preuve.}
On montre cette implication en appliquant tout d'abord 
la conjecture \cite[Appendix]{KR1}
au complexe de Dolbeaut de $f$ et en utilisant le fait que la torsion analytique
relative \'equivariante de ce complexe s'annule.
On calcule ensuite la transform\'ee de Fourier (pour l'action de $\Gal(K/\mQ)$)
de la formule obtenue, puis en utilisant la formule des traces relatives \`a valeurs dans
la cohomologie, on exprime la partie d\'ependant des points
fixes en termes de traces sur la cohomologie de Dolbeaut de $f$. 
Un argument combinatoire permet alors de conclure.
\vskip 0.1cm

Lorsque $n=1$, comme conséquence immédiate du théorème \ref{res1}, 
ou plus simplement grâce à un cas particulier de 
la formule de Lefschetz conjecturale \cite[Appendix]{KR1}
d\'emontr\'e dans \cite{KR}, on peut affirmer que~:
\begin{cor}
Sous les hypoth\`eses du théorème \ref{res1}, 
la conjecture \ref{conjIIa}
vaut apr\`es produit tensoriel par $\mQ(\chi,\mu_q)$ 
lorsque $n=1$. 
\label{corII}
\end{cor}
{\bf Remarque.} Cet énoncé \'etablit en particulier 
une variante de la conjecture \og des p\'eriodes\fg\ de Gross-Deligne  
pour une certaine classe de structures de 
Hodge CM d\'ecoup\'ees dans la cohomologie de vari\'et\'es 
sur $\Qb$.
\begin{theor} La conjecture \cite[Appendix]{KR1}
implique que la conjecture \ref{conjI} 
vaut, apr\`es produit tensoriel par $\mQ(\chi,\mu_{f_K})$, 
pour tout $n \geqslant 1$ lorsque
$S$ est vide et $K$ est une extension ab\'elienne de $\mQ$ 
dont le conducteur $f_K$ est inversible dans $R$. 
\label{res2}
\end{theor} 
{\bf Esquisse de la preuve.} On utilise la th\'eorie 
du corps de classes pour se ramener \`a
$K=\Q(\mu_{f_K})$ puis l'on montre que l'on peut d\'eduire 
ce cas d'une variation de \ref{res1}.
\vskip 0.1cm

Lorsque $n=1$, comme conséquence immédiate du théorème \ref{res2}, 
ou plus simplement grâce au cas particulier de \cite[Appendix]{KR1} 
d\'emontr\'e dans \cite{KR}, on peut affirmer que~:
\begin{cor}
La conjecture \ref{conjI} est vraie apr\`es produit tensoriel par 
$\mQ(\chi,\mu_{f_K})$ lorsque
$n=1$, $S$ est vide et $K$ est une extension ab\'elienne de $\mQ$ 
dont le conducteur $f_K$ est inversible dans $R$.
\end{cor}

\medskip

{\bf Remerciements.} Nous remercions chaleureusement J.-I. Burgos 
et U. K\"uhn pour avoir bien voulu nous expliquer 
leur th\'eorie de l'intersection arithm\'etique avec singularit\'es 
logarithmiques (voir aussi \cite{BuLe} et 
\cite{KuLe}). Leurs commentaires ainsi que ceux   
de J. Bellaiche, P. Colmez, D. Harari et C. Soul\'e nous ont 
\'et\'e fort utiles. Nos remerciements vont \'egalement \`a K. K\"ohler
pour nous avoir communiqu\'e \cite{K1} et pour avoir
attir\'e notre attention sur un point de calcul qui nous avait \'echapp\'e.


\begin{thebibliography}{13}

\bibitem{BiMa}
Bismut, J.-M., Ma, X.:
Equivariant families of immersions and 
analytic torsion forms. A para\^\i tre.

\bibitem{Bo}
Bost, J.-B.:
Potential theory and Lefschetz theorems for arithmetic surfaces.
Ann. Sci. \'Ecole Norm. Sup. (4) {\bf 32}, no. 2, 241--312 (1999).

\bibitem{BuLe}
Burgos, J.-I.:
Lettre aux auteurs (mars 2001).

\bibitem{Bu2}
Burgos, J.-I.: 
Arithmetic Chow rings and Deligne-Beilinson cohomology. 
J. Algebraic Geom. {\bf 6}, 335--377 (1997).

\bibitem{BuWa}
Burgos, J.-I., Wang, S.:
Higher Bott-Chern forms and Beilinson's regulator.
Invent. Math. {\bf 132}, no. 2, 261--305 (1998).

\bibitem{C}
Colmez, P.: P\'eriodes des vari\'et\'es ab\'eliennes \`a multiplication
complexe.
Ann. of Math. (2) {\bf 138}, no. 3, 625--683 (1993).

\bibitem{DeMu}
Deninger, C., Murre, J.: 
Motivic decompositions of abelian schemes and the Fourier 
transform. J. Reine Angew. Math. {\bf 422}, 201--219 (1991). 
 
\bibitem{GS2}
Gillet, H., Soul\'e, C.: Arithmetic intersection theory.
Publications Math. IHES {\bf 72} (1990).

\bibitem{GS3}
Gillet, H., Soul\'e, C.: Characteristic classes for algebraic vector
bundles with hermitian
metrics I, II. Annals of Math. {\bf 131}, 163--203, 205--238 (1990).

\bibitem{GriTu}
Griffiths, P., Tu, L.:
Variation of Hodge structure.
Topics in transcendental algebraic
geometry (Princeton, 1981-1982), Ann. of Math. Stud. {\bf 106}, 3--28:
Princeton Univ. Press 1984.

\bibitem{Gross}
Gross, B. H.:
On the periods of abelian integrals and a formula of
Chowla and Selberg. Invent. Math. {\bf 45}, 193--211 (1978).

\bibitem{K1}
K\"ohler, K.:
A Hirzebruch proportionality principle in Arakelov geometry.
Pr\'epublication de l'Institut de math\'ematiques de Jussieu
no. {\bf 284} (avril 2001).

\bibitem{KR}
K\"ohler, K., Roessler, D.: 
A fixed point formula of Lefschetz type in Arakelov geometry I: 
statement and proof. Invent. Math. {\bf 145}, 333--396 (2001). 

\bibitem{KR1}
K\"ohler, K., Roessler, D.:
A fixed point formula of Lefschetz type in Arakelov geometry II:
a residue formula. Ann. Inst. Fourier {\bf 52}, no. 1, 1--23 (2002). 
 
\bibitem{Kudla}
Kudla, S.:
Integrals of Borcherds forms. Pr\'epublication.

\bibitem{KRY}
Kudla, S., Rapoport, M., Yang, T.:
Derivatives of Eisenstein series and Faltings heights. 
Pr\'epublication.  

\bibitem{Ku}
K\"uhn, U.: Generalized arithmetic intersection numbers. 
J. Reine angew. Math. {\bf 534}, 209--236 (2001).

\bibitem{KuLe}
K\"uhn, U.: Letter to the authors (avril 2001).

\bibitem{KuThe}
K\"uhn, U.: \"Uber die arithmetischen Selbstschnittzahlen zu Modulkurven
und Hilbertschen Modulfl\"achen. Dissertation Hu-Berlin (1999).

\bibitem{Mu}
Mumford, D.:
Hirzebruch's proportionality theorem in the noncompact case.
Invent. Math. {\bf 42}, 239--272 (1977).

\end{thebibliography}
\end{document}